\documentclass[12pt]{article}
\usepackage{amssymb}
\usepackage{amsfonts}
\usepackage{amsmath}
\usepackage{amsthm}

\title{Classical versus quantum probability: Comments on the paper ``On universality of classical probability with contextually labeled random variables''\\ by E. Dzhafarov and M. Kon}

\author{Andrei Khrennikov\\International Center for Mathematical Modeling \\
in Physics, Engineering, Economics, and Cognitive Science\\
Linnaeus University, V\"axj\"o, Sweden  }

\begin{document}
\maketitle

\begin{abstract} Recently Dzhafarov and Kon published the paper advertising  the possibility to use the coupling
 technique of classical probability theory 
to model incompatible observables in quantum physics and quantum-like models of psychology.
 Here I present comments on this paper by stressing advantages and disadvantages. 
\end{abstract}

\section{Introduction}

In a series of papers  Dzhafarov and coauthors present applications of the {\it coupling technique of classical probability} (CP) theory to modeling incompatible observables 
in quantum physics  and psychology
(see references in Dzhafarov \& Kon, 2018; Dzhafarov \& Kujala, 2016; Dzhafarov, Kujala, Cervantes, Zhang, \& Jones, 2016;
 Dzhafarov, Kujala, \& Larsson, 2015;  Dzhafarov, Zhang, \& Kujala, 2015).  Coupling  
is a well developed CP-technique (Lindvall,1992;  Thorisson, 2000) for unifying generally unrelated random variables 
on the basis of a common 
Kolmogorov probability space (Kolmogorov, 1952).  We recall that in  CP-applications, observables are
 represented by random variables. The crucial contribution of Dzhafarov and coauthors
is the proposal to connect the coupling CP-technique with incompatible observables of quantum physics as well as 
psychology. Dzhafarov and coauthors called 
their approach {\it Contextuality-by-Default} (CbD).  

 We recall that observables are incompatible 
if they cannot be measured jointly. In the probabilistic terms it means that  their {\it  joint probability distribution} (jpd) 
does not exist. Instead of the jpd, one has to operate 
with a family of probability distributions depending on experimental contexts 
(see the  {\it V\"axj\"o model} in Khrennikov, 2010). This is a good place 
to make a remark about terminology. We do not consider the V\"axj\"o model  to be a part of CP. Here we operate 
with a family of Kolmogorov probability spaces labeled 
by contexts (Kolmogorov,  1956). These probability spaces are consistently coupled with the aid of
 transition probabilities. In the simplest case of dichotomous observables
it is possible to present conditions on  transition probabilities implying the possibility 
to represent the contextual model in the complex Hilbert space (Khrennikov, 2010).
 
The main advantage of Dzhafarov and coauthors'  approach to contextuality is its applicability 
 to statistical data violating the condition known in quantum 
physics as {\it no-signaling} (or no-disturbance, parameter invariance,  etc. --- in CbD it is 
called ``consistent connectedness'').
 The physical terminology, no-signaling,  is quite ambiguous.  Mathematically  no-signaling means consistency of 
marginal probabilities for different and generally incompatible contexts. 

We also point to a different approach to  CP-unification of observables related to generally incompatible contexts,
proposed  by Khrennikov (2015). It is based on the standard procedure of randomization. 

Personally,  I am sympathetic to modeling contextual statistics of generally incompatible measurements by using a 
classical measure-theoretic approach,  and I contributed a good deal to such studies (see, e.g.,  Khrennikov, 2010, 2015).  
Application of the coupling method is an important step in this direction, especially towards demystification 
of the use of  {\it quantum probability} (QP) in physics or in modeling cognitive, social and political processes.
However, I think that QP is a better mathematical  tool for modeling these processes.\footnote{For examples of use of QP in psychology, see  
monographs  Asano,  Khrennikov, 
 Ohya,   Tanaka,   and Yamato (2015),  Bagarello (2012). Busemeyer and  Bruza (2012), 
Haven and Khrennikov (2013),  Khrennikov (2004a, b, 2010);  handbook  Haven and Khrennikov (2017); 
 textbook Haven, Khrennikov, and Robinson (2017) and some other 
representative papers  (Aerts, 2009; Asano, Masanori, Tanaka, Basieva,   \& Khrennikov, 2011;  
 Bagarello,   Basieva,  Pothos, \&   Khrennikov, 2018;
 Basieva,  Pothos,  Trueblood,  Khrennikov,  \& Busemeyer, 2017;   
Busemeyer,  Wang, \& Townsend, 2006; 
 Busemeyer,  Wang, \&  Mogiliansky-Lambert, 2009;   Haven \&  Khrennikov, 2009; Haven  \& Sozzo,  2016; 
 Khrennikov, 1999; Khrennikova, 2014, 2017;   Pothos \&  Busemeyer, 2009,  2013).}  

\section{Complementarity and the context-dependence of probability}

{\bf Bohr's  principle of complementarity.}    
An output of  any observable   is composed of the contributions of a system and a measurement device.
The whole experimental context has to be taken into account.
There is no reason to expect that all experimental contexts can be combined. 
Therefore one cannot expect that all observables can be measured jointly.
There exist incompatible observables.   See Plotnitsky (2012) for details.

{\bf Complementarity as contextuality of probability.}  The  principle of complementarity can be reformulated in probabilistic terms. In short, we can say that 
{\it the measurement part of QM is a calculus of context-dependent probabilities.} This viewpoint was presented
 in a series of works by the author of this comment (e.g., Khrennikov, 2010).
We emphasize  that QP  is a very special contextual probabilistic calculus. {\it Its specialty 
consists in using  the wave function $\psi$ to unify
generally incompatible contexts.}

In classical statistical physics the contextuality of observations is not emphasized. Here it is assumed 
that it is possible to introduce a single 
context-independent probability measure $P$ and reproduce the probability distributions of all physical 
observables on the basis of $P.$\footnote{ Mathematically the observables are presented
by random variables that are functions on the space of {\it elementary events} $\Omega.$  
Events (representing the outputs of observations) 
are represented as subsets of $\Omega.$ The set of events  ${\cal F}$ is endowed with the
 structure of a set $\sigma$-algebra, i.e., it is closed with respect to the operations of 
conjunction, disjunction, and negation. These operations are represented as intersection,
 union, and complement  (Kolmogorov, 1956, but see already Boole, 1862, 1958).}
 
{\bf Non-existence  of the joint probability distribution.}
Let  ${\cal P}= (\Omega, {\cal F}, P)$  be a Kolmogorov probability space (Kolmogorov, 1956). Each random variable $a: \Omega \to \mathbf{R}$ determines 
the probability distribution $P_a.$ The crucial point is that all these distributions are encoded in
 the same probability measure $P: P_a (\alpha) = P(\omega \in \Omega: a(\omega) = \alpha).$ Thus, {\it probability distributions of all observables (represented by 
random variables) can be consistently unified on the basis of $P.$}  
For any pair of random variables $a, b,$  their jpd $P_{a, b}$ is defined and the following condition of marginal consistency holds:
\begin{equation}
\label{MC}
P_a(\alpha)= \sum_\beta P_{a,b} (\alpha, \beta)
 \end{equation}
This condition means that observation of $a$  jointly with  $b$ does not change the probability distribution of  $a.$ 
 Equality (\ref{MC}) implies that, for any two observables $b$ and $c,$ 
\begin{equation}
\label{MC1}
\sum_\beta P_{a,b} (\alpha, \beta) = \sum_\gamma P_{a,c} (\alpha, \gamma).
 \end{equation}
 In fact, condition (\ref{MC1}) is equivalent to (\ref{MC}):  by selecting the random variable $c$  such that $c(\omega)=1$ almost everywhere, we see that (\ref{MC1}) implies (\ref{MC}). 
These considerations are easily generalized to a system of $k$ random variables $a_1,..., a_k.$ Their jpd  is well defined,
$P_{a_1,...,a_k}(\alpha_1,..., \alpha_k)= P(\omega \in \Omega: a_1(\omega) = \alpha_1, ...., a_k(\omega) = \alpha_k).$     

Consider now some system of {\it experimental observables} $a_1,..., a_k.$ If  the experimental design 
for their joint measurement exists, then 
it is possible to define their jpd $P_{ a_1,...,  a_k}(\alpha_1,..., \alpha_k)$ (as the relative frequency of their joint outcomes). 
This probability measure $P \equiv P_{a_1,...,  a_k}$
can be used to define the Kolmogorov probability space, i.e., the case of joint measurement can be described by CP. 

Now consider the general situation: only some groups of observables can be jointly measured. For example, there are four observables $a_1,  a_2$ and $ b_1,  b_2$
and we are able to design  measurement procedures only  for some pairs of them, say $(a_i, b_j), i, j=1,2.$ In this situation, there is no reason to expect
that one can define mathematically the joint probability
distribution $P_{a_1, a_2, b_1, b_2}(\alpha_1, \alpha_2, \beta_1, \beta_2)$ such that the conditions of the marginal consistency for pairs hold:
\begin{equation}
\label{MCx}
P_{a_1, b_1}(\alpha_1, \beta_1) = \sum_{\alpha_2, \beta_2} P_{a_1, a_2, b_1, b_2}(\alpha_1, \alpha_2, \beta_1, \beta_2),....
\end{equation}
This situation is typical for quantum theory.   This  is a complex interplay 
of theory and experiment. Only probability distributions $P_{a_i, b_j}$ can be experimentally verified. The jpd 
$P_{a_1, a_2, b_1, b_2}$ is a {\it hypothetical mathematical quantity.} However, if it existed, one may expect that 
there would be  some 
experimental design for joint measurement of the quadruple of observables $( a_1, a_2, b_1, b_2).$ On the other hand, 
if it does not exist, then it is meaningless even to try to design an experiment for their joint measurement.  

{\bf The CHSH inequality.}
How can one get to know whether a jpd exists?  The answer to this question is given by a theorem (Fine, 1982) 
concerning one of the Bell-type inequalities 
(Bell, 1964, 1987), namely,  the CHSH inequality  (Clauser, Horne, Shimony, \& Holt , 1969).
Consider  covariation of  compatible observables $a_i$ and $b_j$ given by  $\langle  a_i  b_j\rangle = E [a_i b_j]= \int  \alpha \beta \;  dP_{a_i, b_j}(\alpha, \beta).$
By  Fine's theorem a jpd exists if and only if the CHSH-inequality for these correlations is satisfied, namely, the inequality
\begin{equation}
\label{CH}
\vert \langle  a_1  b_1\rangle  - \langle  a_1  b_2\rangle + \langle  a_2  b_1\rangle +\langle  a_2 b_2\rangle   \vert \leq 2.
\end{equation}
and the three other inequalities corresponding to all possible permutations of indexes $i=1,2$ and $j=1,2.$  

We restrict further  considerations to the CHSH-framework, i.e., we shall not consider other types of Bell inequalities.  

 The above presentation of Fine's result  is common for physics' folklore. However,  Fine did not consider explicitly  the CHSH inequalities presented above, see (\ref{CH}). 
He introduced four inequalities that are necessary and sufficient for 
the jpd to exist, but these inequalities are expressed differently to the CHSH inequalities. 
The  CHSH inequalities are derivable from Fine's four inequalities stated in Theorem 3 of his paper.

{\bf Remark.} In quantum physics this very clear and simple meaning of violation of the CHSH-inequality  
is obscured by the issue of  nonlocality.
Non-locality is relevant to space separated systems. However, except for perhaps cognitive neuro-science, 
cognitive psychology does not model space separated systems. 
 We remark  that Bell type inequalities were considered already by Boole (1862, 1958) as necessary 
conditions for existence of a jpd.  

The above reasoning has an important consequence. The existence of a jpd implies that the condition 
of marginal consistency with respect to the jpd  should hold not only for the pairwise 
probability distributions (see (\ref{MC})), but even for probability distributions of each observable, $a_i, b_j, i,j =1,2,$ i.e., 
\begin{equation}
\label{MCT}
P_{a_i}(\alpha) = \sum_{\beta} P_{a_1, b_j}(\alpha,  \beta), j=1, 2, \;  P_{b_j}(\beta) = \sum_{\alpha} P_{a_i, b_j}(\alpha,  \beta),  i=1,2.
\end{equation}
As was pointed out, this condition  is known in quantum physics as  the {\it no-signaling condition.} 
 Thus, the Fine theorem presupposed that two conditions of marginal consistency, (\ref{MCx}) and  (\ref{MCT}), 
 jointly hold. 

{\bf Signaling in physical  and  psychological experiments.}
By using the quantum calculus of probabilities, it is easy to check whether the no-signaling condition holds for quantum observables, which are represented mathematically by Hermitian operators.  
Therefore Fine's theorem is applicable to quantum observables. This theoretical fact 
played an unfortunate role in hiding from view  (\ref{MCT}) in experimental research on violation 
of the CHSH-inequality. Experimenters were focused on observing as high  violation of (\ref{CH}) 
as possible and they ignored the no-signaling condition (\ref{MCT}).
However, if the latter is violated, then a jpd automatically does not exist, and there is no reason to expect 
that  (\ref{CH}) would be satisfied. 
The first paper in which the signaling issue  in  quantum experimental research was highlighted was Adenier and 
Khrennikov (2007). There it was shown that  statistical data collected in 
the basic experiments (for that time) performed by Aspect (1983)  and Weihs (1999) violates the no-signaling condition. 
\footnote{After this publication experimenters became aware of the signaling issue and started 
to check it (Giustina et al.,  2015, Shalm et al., 2015). However, analysis presented in Adenier and 
Khrennikov (2016)  demonstrated that even statistical data 
generated in the first loophole-free experiment to violate the CHSH-inequality (Hensen et al., 2015) 
exhibits very strong signaling. In Section \ref{MO}  the discussion on signaling in physics and psychology will be continued.
}

The experiments to check CHSH and other Bell-type inequalities have also been performed for mental observables in the form of questions asked to people. The first such experiment 
was done by Conte, Khrennikov, Todarello, Federici, Mendolicchio, and  Zbilut (2008) 
and was based on the theoretical paper Khrennikov (2004a).
 As was  found by Dzhafarov,  Zhang, and Kujala  (2015) all known experiments 
of this type suffer of signaling. Moreover, in contrast to physics, in psychology there are no 
theoretical reasons to expect no-signaling. In this situation Fine's theorem is not applicable. 
And  Dzhafarov and his coauthors were the first who understood  the need of adapting 
 the Bell-type inequalities to experimental data exhibiting signaling.
Obviously, the interplay of whether or not a jpd exists for quadruple 
$$
S=(a_1, a_2, b_1, b_2)
$$  
can't be considered  for signaling data.

{\bf Coupling method (Contextuality-by-Default).} Dzhafarov and his coauthors (see references in the introduction)  propose considering, instead of quadruple $S,$ some octuple ${\bf S}$
generated by doubling each  observable and associating ${\bf S}$ with four contexts of measurements of pairs,
 $C_{11}=(a_1, b_1), C_{12}=(a_1, b_2), C_{21}=(a_2, b_1), C_{22}=(a_2, b_2).$ 
Thus, the basic object of CbD-theory has the form 
$$
{\bf S}= (A_{11}, B_{11},  A_{12}, B_{21}, A_{21}, B_{12}, A_{22}, B_{22}).
$$
It is assumed that this system of observables can be realized by random variables on {\it the same Kolmogorov 
probability space} ${\cal P}_{{\bf S}}= ({\bf \Omega}, {\cal \bf{F}},  \bf P).$
(We shall use bold symbols for sample spaces and probabilities realizing the octuple representation 
of observables by random variables.) 
For example,  $A_{ij}= A_{ij}(\omega), \omega \in\bf  \Omega,$ is a  random variable representing 
observable $a_i$ measured jointly with  the observable $b_j.$  

By moving from quadruple $S$  to octuple ${\bf S},$ one confronts the problem of identity of an observable
which is now represented by two different random variables, e.g., 
the observable $a_i$ is represented by  the random variables $A_{ij}(\omega) , j=1,2.$
In the presence of signaling one cannot expect the equality of two such random variables almost everywhere. 
Dzhafarov and coauthors (Dzhafarov \& Kujala, 2016; Dzhafarov, Kujala, Cervantes, Zhang, \& Jones, 2016;
 Dzhafarov, Kujala, \& Larsson, 2015;  Dzhafarov, Zhang, \& Kujala, 2015)  
came  up with a novel treatment of the observable-identity problem.  

It is assumed that averages
\begin{equation}
\label{Ct}
m_{a; ij}=\langle A_{ij} \rangle, \; m_{b; ij}= \langle B_{ij}\rangle
\end{equation}
and covariation
\begin{equation}
\label{Ct1}
{\cal C}_{ij} =\langle A_{ij}  B_{ji} \rangle 
\end{equation}
are fixed. These are measurable quantities.  They can be statistically verified by experiment. 

Set 
\begin{equation}
\label{Cty}
\delta(a_i)= m_{a; i1}- m_{a; i2} \; \delta(b_j)= m_{b; j1} - m_{b; j2}, 
\end{equation}
and 
\begin{equation}
\label{Cty1}
\Delta_0 = \frac{1}{2} \Big(\sum_i \delta(a_i) + \sum_j \delta(b_j)\Big).
 \end{equation}
This is the experimentally verifiable measure of signaling.  

We remark that in the coupling representation the joint satisfaction of the CHSH inequalities, i.e.,  (\ref{CH}) and other
inequalities obtained from it via permutations, can be written in the form:
\begin{equation}
\label{CHBJ}
\max_{ij} \vert \langle  A_{11}  B_{11}\rangle +
\langle  A_{12}  B_{21}\rangle + 
\langle  A_{21}  B_{12}\rangle +
\langle  A_{21}  B_{22}\rangle - 2 \langle  A_{ij}  B_{ji} \rangle\vert  \leq 2.
\end{equation}
In the signaling-free situation, e.g., in quantum physics,  the difference between the left-hand and right-hand sides is considered as 
the measure of contextuality. Denote (1/2 times) this quantity by $\Delta_{\rm{CHSH}}.$ It is also experimentally verifiable.     

Then Dzhafarov and coauthors introduced quantity
\begin{equation}
\label{CHPP0}
\Delta({\bf P})= \sum \Delta_{a_i}({\bf  P})+ \sum \Delta_{b_j}({\bf P}),
\end{equation}
where 
\begin{equation}
\label{CHPP}
\Delta_{a_i}({\bf P})=  
{\bf  P}(\omega: A_{i1}(\omega) \not = A_{i2}(\omega) ), \Delta_{b_j}({\bf P})= 
{\bf P}(\omega: B_{j1}(\omega) \not  = B_{j2}(\omega)).
\end{equation}
Here  $\Delta_{a_i}({\bf P})$ characterizes mismatching of representations of observable $a_i$ 
by  random variables  $A_{i1}$ and $A_{i2}$ with respect to probability measure ${\bf P};$
   $\Delta_{b_j}({\bf  P})$ is interpreted in the same way.
The problem of the identity of  observables is formulated as the mismatching minimization or identity maximization problem   
\begin{equation}
\label{CHPP0}
\Delta({\bf P}) \to \min
\end{equation}
with respect to all octuple probability distributions ${\bf P}$ satisfying constraints (\ref{Ct}), (\ref{Ct1}). 
And it turns out, that 
\begin{equation}
\label{CHPP0yy}
\Delta_{\rm{min}}=\min \Delta({\bf P})= \max [\Delta_0,   \Delta_{\rm{CHSH}}].
\end{equation}
It is natural to consider the solutions of the identity maximization problem (\ref{CHPP0}) as 
CP-representations for contextual system ${\bf S}.$ 
The corresponding random variables have the highest possible, in the presence of signaling,  degree of identity.  

The quantity $\Delta_{\rm{min}} - \Delta_0$ is considered as the measure of ``genuine contextuality''. Thus contextuality encoded in $\Delta_{\rm{CHSH}}$ is the sum of 
``straightforward contextuality'' carried via signaling  and genuine contextuality. This approach is very useful to study contextuality in the presence of signaling.      
The key point is the coupling of this measure of contextuality with the problem of the identity of observables measured in different contexts.  As was emphasized by 
Dzhafarov \& Kujala (2015)\footnote{In my opinion this is the best paper about CbD 
combining clarity and simplicity with rigorousness of presentation.}: 

``contextuality means that random variables recorded under mutually
incompatible conditions cannot be "sewn together" into a single system
of jointly distributed random variables, provided one assumes
that their identity across different conditions changes as little as
possibly allowed by direct cross-influences (equivalently, by observed
deviations from marginal selectivity).''

This approach to contextuality can be reformulated in the CHSH-manner by using what we can call Bell-Dzhafarov-Kujala (BDK) inequality.

 the BDK-inequality
\begin{equation}
\label{CHBJ}
\max_{ij} \vert \langle  A_{11}  B_{11}\rangle +
\langle  A_{12}  B_{21}\rangle + 
\langle  A_{21}  B_{12}\rangle +
\langle  A_{21}  B_{22}\rangle - 2 \langle  A_{ij}  B_{ji} \rangle \vert  - 2\Delta_0  \leq 2.
\end{equation}
It was proven that  octuple-system ${\bf S}$ exhibits no genuine contextuality, i.e., 
\begin{equation}
\label{CHBJz}
\Delta_{\rm{min}}=\Delta_0,
\end{equation}
  if and only if the BDK-inequality is satisfied.
  
 {\bf Non-uniqueness.}  
We point out the following  feature of the CbD-modeling 
 physical and psychological phenomena:  {\it a coupling guaranteeing maximal possible  
identification of different classical random variables representing the same observable is not unique.}  
Optimization problem (\ref{CHPP0}) has a non-unique solution with a high degree of degeneration.

{\it How can one select the ``right coupling'? }

It is seems that the ``right coupling'' does not exist. By using the CbD-approach we go beyond observational theories such as, for example,  quantum theory.
Consider probability space ${\cal P}_{{\bf S}}= ({\bf \Omega}, {\cal \bf{F}},  \bf P)$ for octuple ${\bf S}.$ Elementary events of this space, $\omega \in {\bf \Omega},$
are ``hidden variables''. These hidden variables are contextual. It is well known that even in the absence of signaling, e.g., 
in quantum theory, plenty of contextual hidden variable models matching observational data can be constructed.
 One explanation of this multitude of models is that real physical (or psychological)
contextuality is determined not only by semantically defined observables, but also by apparatuses  used for their measurement. The same observable determining 
context can be measured by a variety of apparatuses. These context-apparatuses are represented by probability spaces generated in the CbD-approach.

In contrast to the CbD-theory,  {\it QP does not suffer the non-uniqueness problem. There is one fixed  quantum state given by a normalized vector $\psi$ or generally 
by a density operator $\rho,$ and there is the unique representation of observables by Hermitian operators.}
The QP-description is the natural generalization of the CP-description based on the single probability measure $P.$ 
 In particular, by applying QP to cognition and psychology we can identify 
quantum states with mental states and obtain a consistent model of decision making based on such quantum-like states. 

\section{Classical, contextual, and quantum probability}

The main message of Dzhafarov and Kon (2018) is that CP can be used to describe mathematically statistical data collected in all possible 
experiments in physics and psychology.  And they rightly pointed out that some authors have actively claimed that generally 
CP is inapplicable to some experimental data (see references in  Dzhafarov \& Kon, 2018).  
I generally support this critique and agree that often statements about inapplicability of CP were formulated vaguely. 
At the same time some  authors mentioned in Dzhafarov and Kon (2018),  including 
Feynman and myself,  understood the interrelation between CP and QP
very well.

To clarify  their position, one has to recall that by 
the Copenhagen interpretation quantum mechanics is an observational theory and all its statements have to be formulated for 
experimentally verifiable data.   The corresponding statements about inapplicability of CP are 
about the impossibility of defining a Kolmogorov probability space which is based solely on experimentally  verifiable events. This is the viewpoint of Feynman  who pointed  out 
that (in the observational framework) one can describe the two-slit experiment only by violating the additivity of probability. He stressed 
that this experiment is, in fact, composed of three different experiments: both slits are open and two experiments in which only 
one  of the slits is open. Denote these experimental contexts $C_{12}, C_{1}, C_{2}.$ For each of these contexts, CP works very well.
However,  it is impossible to represent observables with respect to these contexts in the same Kolmogorov probability space without 
introducing unobservable events. Then Feynman proposed to employ a generalized probabilistic model with non-additive probabilities. 
The author of this comment reformulated Feynman's considerations by using conditional probabilities. This approach led to a modification 
of the formula of total probability: perturbation of the CP-formula by an interference term (Khrennikov, 2004a, b, 2010).   

By moving from the purely observational theory (such as quantum mechanics) to, so to say, sub-observational models, we can proceed with CP (see, e.g.,  Khrennikov, 2014) 
for so-called {\it prequantum classical statistical field theory} (PCSFT). 

The main problem is that the majority of scientists do not separate the two layers of mathematical modeling of natural and mental phenomena: theoretical
 and observational (see  Atmanspacher \&  Primas,  2005;  Bolzmann, 1905; Hertz, 1899;   
  Khrennikov, 2017, 2018).  Theoretical and observational models for some phenomena 
are coupled via some correspondence mapping for states and variables.  The CbD and PCSFT models are theoretical. 
 Feynman and I  discussed impossibility to proceed with CP at the observational level.

Regarding statements of other authors mentioned in  Dzhafarov and Kon (2018) and claiming inapplicability of CP in some situations, it is difficult to say whether their 
statements were about theoretical or observational models or their mixing. The latter is the most common, 
since people typically do not explicitly identify  the  modeling layer.   
 
 \section{Impact of CbD theory}
 
 The main impact of the CbD theory is the possibility of introducing a measure of contextuality for statistical data with signaling. Generalization of no-signaling contextuality 
is done very smoothly with the BDK-inequality playing the role of the CHSH-inequality.\footnote{Note that we restrict our considerations to the measurement scheme related to the 
CHSH-inequality. Generally the CbD theory is applicable to all known measurement schemes for discrete observables.}

We remark that the CbD approach, i.e., the use of the coupling technique of CP, has a nontrivial impact even in the absence of signaling.  It demystifies quantum mechanics by highlighting 
the  role of contextuality, i.e., dependence of observables on the whole ``experimental arrangement'', as was  emphasized by Bohr (see also Khrennikov, 2004a, b, 2010). 

Finally, I would like to point to one very important consequence of the possibility of the CP-description of complexes of contexts, such as in the two-slit or 
 Bell-type experiments. Experimental data is the subject of the statistical analysis. The latter is fundamentally based on CP, in particular, the representation  of observables
by random variables. Thus, to justify mathematically the statistical significance of a
violation the CHSH-inequality, one has to proceed on the basis of some Kolmogorov probability 
space. The CbD theory provides such a possibility. 

\section{Mental signaling: fundamental or technical?}
\label{MO}

Dzhafarov, Zhang,  and Kujala, (2015) pointed out  that statistical data collected in psychological 
experiments contains statistically significant 
signaling patterns. One can wonder whether signaling is a fundamental feature of mental observations or  a mere 
technicality, perhaps the consequence of badly designed 
or/and performed experiments. We recall that in physics signaling patterns 
were found in all Bell experiments during  the first 30 years. 

Since quantum theory predicts the absence of signaling, signaling patterns in experimental statistical 
data  were considered to be a technicality.\footnote{Here ``technicality'' refers to situations in which technical
equipment, experimental design, improper calibration of detectors and so on,
influence an experiment's results.} 
Understanding the technical sources of signaling and finding ways  to eliminate it 
 required great efforts of the experimenters. 
Finally, Giustina et al. (2015) and Shalm et al. (2015)  reported that the null hypothesis of signaling can
 be rejected for the data collected in these
experiments.\footnote{Hensen et al.  (2015) also claimed that signaling hypothesis can be rejected. 
However, the independent analysis of their 
statistical data performed by Adenier and Khrennikov (2016) showed the presence of statistically significant 
signaling. It is very important to perform a similar independent analysis for the 
data obtained by Giustina et al. (2015) and Shalm et al. (2015).   }  

In psychology  the situation is more complicated. 
There are no theoretical reasons to expect no signaling. Therefore, it is not obvious whether
signaling is a technicality or a fundamental   feature of cognition. For the moment, only a few experiments have been 
performed. One cannot exclude that in the future  more advanced 
experiments would generate data without signaling. As the first step towards such experiments, 
possible experimental sources of mental signaling should be analyzed.  
However,   it may be that mental signaling is really a fundamental feature of cognition.

In any event,  it is interesting to attempt  to find non-signaling contextual patterns in human behavior. If such patterns were found, then one can try to connect 
such non-signaling contextuality with some specialty in human psychology. 

\section*{Acknowledgments}

I would like to that thank E. Dzhafarov for his numerous comments which were helpful to understand correctly his position, I also thankful to I. Basieva and E. Pothos for 
long and exciting discussions on quantum versus psychological foundations during my sabbatical  at City, University of London. 
         
{\bf References}

Adenier, G., \&  Khrennikov, A.  (2007).  Is the fair sampling assumption supported by EPR experiments?  
{\it Journal of  Physics B: Atomic, Molecular and Optical Physics,  40}, 131-141.

 Adenier, G.,  \&  Khrennikov, A. (2016).  Test of the no-signaling principle in the Hensen loophole-free CHSH experiment. 
{\it Fortschritte der Physik (Progress in Physics),  65},  1600096.

Aerts, D. (2009). Quantum structure in cognition. \textit{Journal of Mathematical Psychology, 53}, 314-348.

Asano, M., Masanori. O., Tanaka, Y., Basieva, I.,   \& Khrennikov, A. (2011). Quantum-like model of brain's functioning: 
Decision making from decoherence {\it Journal of Theoretical Biology, 281}, 56-64. 

Asano, M., Khrennikov, A.,   Ohya, O.,  Tanaka, Y.,  \& Yamato,  I. (2015). 
 \textit{Quantum adaptivity in biology: from genetics to cognition.}  Heidelberg-Berlin-New York: Springer.

Aspect,  A. (1983). {\it Three experimental tests of Bell inequalities by 
the measurement of polarization correlations between photons.} Orsay: Orsay Press.

Atmanspacher, H.,  \&  Primas, H.  (2005).  Epistemic and ontic quantum realities. 
In  G. Adenier   \&  A. Yu. Khrennikov  (Eds.),  {\it Foundations
of Probability and Physics-3} 750 (pp. 49-62).   Melville, NY: AIP.

Bagarello,  F.  (2012). {\it Quantum dynamics for classical systems: with applications of the number operator.} 
New York: Wiley.

Bagarello, F.,  Basieva, I.,  Pothos, E.,  \&  Khrennikov, A.   (2018). 
Quantum like modeling of decision making: Quantifying uncertainty with the aid of Heisenberg-Robertson inequality.
{\it Journal of Mathematical Psychology, 84}, 49-56.

Basieva, I.,  Pothos, E., Trueblood, J., Khrennikov, A., \& Busemeyer, J.  (2017).  Quantum probability updating 
from zero prior (by-passing Cromwell's rule).  
{\it  Journal of Mathematical Psychology, 77}, 58-69. 

Bell, J.  (1964). On the Einstein-Podolsky-Rosen paradox. {\it Physics, 1}, 195-200.

 Bell J. (1987) {\it Speakable and unspeakable in quantum mechanics.} Cambridge: Cambridge Univ. Press.

Boltzmann, L. (1905). \"Uber die Frage nach der objektiven Existenz der Vorg\"ange in 
der unbelebten Natur. In J. A. Barth (Ed.),   Leipzig: Popul\"are Schriften.   

Boole, G.  (1862). On the theory of probabilities. {\it Philosophical Transactions of the Royal Society of London, 152}, 225-242.

Boole, G.  (1958). {\it An Investigation of the laws of thought.} New York: Dover.

Busemeyer, J. R., Wang, Z., \& Townsend, J. T. (2006).  Quantum dynamics of human decision making.
\emph{Journal of Mathematical Psychology, 50}, 220-241.

 Busemeyer,  J. R., \&  Bruza, P. D. (2012). {\it Quantum models of cognition and decision.} 
Cambridge: Cambridge Press.
 
 Busemeyer, J. R., Wang, Z., \&  Lambert-Mogiliansky, A. (2009). Comparison of quantum and Markov models of decision making.
\textit{Journal of Mathematical Psychology, 53}, 423-433.

Clauser, J. F.,  Horne, M. A.,  Shimony, A., \&  Holt,  R. A.  (1969).  
Proposed experiment to test local hidden-variable theories. {\it Physics Review Letters, 23} (15), 880-884.
  
Conte, E., Khrennikov, A., Todarello, O., Federici. A., Mendolicchio, L.,  \&  Zbilut, J. P. (2008). 
A preliminary experimental verification on the possibility of Bell inequality violation in mental states. 
{\it  NeuroQuantology,  6}, 214-221.

Dzhafarov, E. N., \& Kujala, J. V. (2012). Selectivity in probabilistic causality: Where psychology runs into quantum physics.
{\it Journal of Mathematical Psychology, 56}, 54-63.

Dzhafarov, E. N.,  Kujala, J. V., \& Larsson, J-A. (2015). 
Contextuality in three types of quantum-mechanical systems. {\it Foundations of  Physics, 7}, 762-782.

Dzhafarov, E. N., Zhang, R., \& Kujala, J.V. (2015). Is there contextuality in behavioral and social systems?
{\it Philosophical Transactions of the Royal Society: A, 374},  20150099.

Dzhafarov, E. N.,   \&  Kujala, J. V. (2015)  Probabilistic contextuality in EPR/Bohm-type systems with signaling allowed.
In E.  Dzhafarov,  S.  Jordan,  R.  Zhang,  
 \&  V. Cervantes  (Eds.), {\it  Contextuality from Quantum Physics to Psychology}   
(pp. 287-308).  New Jersey: World Scientific Publishing. 

Dzhafarov, E. N., \& Kujala, J. V. (2016). Context-content systems of random variables: The contextuality-by default
theory. {\it Journal of Mathematical Psychology, 74}, 11-33.

Dzhafarov, E. N., Kujala, J. V., Cervantes, V. H., Zhang, R., \& Jones, M. (2016). 
On contextuality in behavioral data. {\it Philosophical Transactions of the Royal Society: A, 374},  20150234.

Dzhafarov, E. N., \&  Kon, M. (2018).  
On universality of classical probability with contextually labeled random variables. {\it  Journal of Mathematical Psychology,
 85}, 17-24.

Fine, A.  (1982). Joint distributions, quantum correlations, and commuting observables. {\it Journal of Mathematical
 Physics, 23}, 1306.

Giustina, M., Versteegh,   M. A. M., Wengerowsky, S., Handsteiner, J., 
Hochrainer, A., Phelan, K., Steinlechner, F.,  
Kofler, J., Larsson, J.-A., Abellan, C., Amaya, W., 
Pruneri, V.,  Mitchell, M. W., Beyer, J., Gerrits, T., Lita, A.E., Shalm,  L.K., Woo Nam, S., 
Scheidl, T., Ursin, R., Wittmann, B.,   \& Zeilinger, A. (2015). 
A significant-loophole-free test of Bell's theorem with entangled photons. {\it Physics Review Letters, 115}, 250401.

Haven, E.,  \&  Khrennikov, A.  (2009).  Quantum mechanics and violation of the sure-thing principle: The use of probability interference and
other concepts. {\it Journal of Mathematical Psychology, 53}, 378-388.

Haven, E., \&  Khrennikov, A. (2013). {\it Quantum social science.} Cambridge: Cambridge University Press.

Haven, E., \& Khrennikov, A. (Eds.) (2017). {\it The Palgrave handbook of quantum models in social science. 
Applications and grand challenges.}   London: Palgrave Macmillan. 

Haven, E., Khrennikov, A., \& Robinson,  T. R.  (2017). 
{\it Quantum methods in social science: A first course.} Singapore: World Scientific Publishing.

Haven, E., \& Sozzo, S. (2016). A generalized probability framework to model economic agents' decisions 
under uncertainty. {\it International Review  Financial Analysis, 47}, 297-303.

Hensen, B.,  Bernien,  H.,  Dreau,  A. E.,  Reiserer,  A., Kalb,  N., Blok,   M. S., Ruitenberg,  J., 
Vermeulen,  R. F. L., Schouten,   R. N., Abellan,   C.,  Amaya,   W.,  Pruneri, V., 
Mitchell,   M. W., Markham,   M., Twitchen, D. J., Elkouss,   D., Wehner,   S.,   Taminiau, T. H.,  \& Hanson, R.  (2015). Experimental loophole-free 
violation of a Bell inequality using entangled electron spins separated by 1.3 km.  {\it Nature, 526}, 682-686.

Hertz, H. (1899). {\it  The principles of mechanics: presented in a new form.} London: Macmillan.

Khrennikov, A.  (1999). Classical and quantum mechanics on information spaces with applications to cognitive, psychological,
social and anomalous phenomena.  {\it Foundations of Physics,  29},  1065-1098.

Khrennikov, A.  (2004a).  On quantum-like probabilistic structure of 
mental information. {\it Open Systems and Information Dynamics,
11}, 267-275.

 Khrennikov, A. (2004b). {\it Information dynamics in cognitive, psychological and anomalous phenomena.} 
Berlin-Heidelberg-New York: Springer.

Khrennikov, A. (2010). {\it Ubiquitous quantum structure: From psychology to finance.} Berlin-Heidelberg-New York: Springer.

Khrennikov, A. (2014). {\it Beyond quantum.}  Pan Stanford Publishing:  Singapore.
  
 Khrennikov, A.   (2015). CHSH inequality: quantum probabilities as classical conditional probabilities. 
{\it Foundations of Physics, 45},   711-725.   

Khrennikov, A. (2017).  Quantum epistemology from subquantum ontology: Quantum mechanics from theory of
classical random fields. {\it Annals of Physics,  377},  147-163.

Khrennikov, A. (2018).	Hertz's viewpoint on quantum theory. 	arXiv:1807.06409 [quant-ph].

Khrennikova, P. (2014).  A quantum framework for `Sour Grapes' in cognitive dissonance. 
\textit{Lecture Notes in Computer Science, 8369}, 270-280.

Khrennikova, P. (2017).  Modeling behavior of decision makers with the 
aid of algebra of qubit creation-annihilation operators.
{\it  Journal of Mathematical Psychology,  78}, 76-85.

Kolmogorov, A. N. (1956). \textit{Foundations of the theory of probability}. New York: Chelsea Publishing Company.

Lindvall, T. (1992).  {\it Lectures on the coupling method.} New York: Wiley.

Plotnitsky, A. (2012). {\it Niels Bohr and complementarity: An introduction.} Berlin and New York: Springer.

Pothos, E. M. , \& Busemeyer, J. R. (2009) A quantum probability explanation for violations of rational decision
making \textit{ Philosophical Transactions of the Royal Society: B, 276},  2171-2178.

Pothos, E. M., \& Busemeyer, J. R. (2013). Can quantum probability provide a new direction for cognitive modeling?
\textit{Behavioral Brain Science, 36}, 255-327.
    
Shalm, L.K., Meyer-Scott, E.,  Christensen, B. G., Bierhorst, P., Wayne,  M. A., Stevens, M. J., 
Gerrits, T., Glancy, S., Hamel, D. R., Allman, M. S., Coakley, K. J.,  
Dyer,  S. D., Hodge, C., Lita,  A. E.,   Verma, V. B.,  Lambrocco, C., Tortorici, E., 
Migdall,  A. L.,  Zhang, Y., Kumor,  D. R., Farr, W. H., Marsili, F., Shaw,  M. D.,  Stern, J. A.,  
Abellan, C., Amaya,  W., Pruneri, W., Jennewein, T., 
Mitchell,  M. W.,   Kwiat, P. G., Bienfang, J. C., Mirin, R. P., Knill, E., \& Woo Nam, S. (2015).  
A strong loophole-free test of local realism.  {\it Physics Review Letters,  115}, 250402.

Thorisson, H. (2000). {\it Coupling, stationarity, and regeneration.} New York: Springer.

Weihs, G.  (1999). {\it Ein Experiment zum Test der Bellschen Ungleichung unter Einsteinscher Lokalit\"at.} 
Vienna: University of Vienna.

\end{document}